\documentclass[12pt]{article}
\usepackage{amssymb,amsthm,amsmath,latexsym}
\usepackage{url}
\newtheorem{theorem}{Theorem}

\newtheorem{proposition}{Proposition}

\theoremstyle{definition}

\newtheorem{remark}{Remark}

\newcommand{\allone}{\mathbf{1}}
\DeclareMathOperator{\rank}{rank}

\begin{document}

\title{The weight distribution of the self-dual $[128,64]$ polarity 
design code\footnote{Research of M. Harada supported by JSPS KAKENHI
Grant Number 15H03633.
Research of E. Novak and V. D. Tonchev supported by NSA Grant 
H98230-15-1-0042}
}

\author{
Masaaki Harada\thanks{
Research Center for Pure and Applied Mathematics,
Graduate School of Information Sciences,
Tohoku University, Sendai 980--8579, Japan.
email: mharada@m.tohoku.ac.jp.},
Ethan Novak\thanks{
Department of Mathematical Sciences,
Michigan Technological University,
Houghton, MI 49931, USA.
email: ewnovak@mtu.edu.}
and 
Vladimir D. Tonchev\thanks{
Department of Mathematical Sciences,
Michigan Technological University,
Houghton, MI 49931, USA.
email: tonchev@mtu.edu.}
}

\maketitle
\date{}

\begin{abstract}

The weight distribution of the binary self-dual $[128,64]$ code
being the extended code $C^{*}$ of the code $C$ spanned by the incidence vectors
of the blocks of the polarity design in $PG(6,2)$ \cite{JT09}
is computed.
It is shown also that $R(3,7)$ and  $C^{*}$
have no self-dual $[128,64,d]$ neighbor with $d \in \{ 20, 24 \}$.

\end{abstract}

\section{Introduction}

We assume familiarity with basic facts and notions from coding theory
and combinatorial design theory (\cite{AK}, \cite{MS},  \cite{T88}).
 
We denote by $PG_{s}(m,q)$ the design
having as points and blocks the points and $s$-subspaces
of the $m$-dimensional projective geometry $PG(m,q)$
over a finite field $GF(q)$ of order $q$,
where $q = p^{t}$ is a prime power and $1 \le s \le m-1$.
The projective geometry design $PG_s(m,q)$ is a 2-$(v,k,\lambda)$
design with parameters
\begin{equation}
\label{pgdmq}
 v = \frac{q^{m+1}-1}{q-1}, \; k = \frac{q^{s+1}-1}{q-1}, \; 
\lambda = \left[ m-1 \atop s-1 \right]_{q}, 
 \end{equation}
where $\left[ m \atop i \right]_{q}$  denotes the  Gaussian coefficient
given by
\[  \left[ m \atop i \right]_{q}=
\frac{(q^m -1) (q^{m-1}-1) \cdots (q^{m - i +1} -1)}{(q^i -1) (q^{i-1}-1)\cdots (q -1)}. \]
The affine geometry design $AG_{s}(m,q)$, ($1\le s\le m-1$), is a
2-$(v,k,\lambda)$ design of the points and $s$-subspaces
of the $m$-dimensional affine geometry $AG(m,q)$ over $GF(q)$, where
\begin{equation}
\label{agdmq}
v = q^{m}, \; k = q^{s}, \; \lambda =  \left[ m-1 \atop s-1 \right]_{q}.
\end{equation}
In the special case when $q=2$, $AG_s(m,2)$, $s\ge 2$, is also a 3-design,
with every three points contained in  $\lambda_3$ blocks, where
\begin{equation}
\label{af2}
 \lambda_3 =\left[ m-2 \atop s-2 \right]_{2}. 
\end{equation}

A finite geometry code (or  geometric code), is a linear code
being the null space of the incidence matrix of a geometric design,
$AG_{s}(m,q)$ or $PG_s(m,q)$. The codes based on affine geometry designs, $AG_s(m,q)$,
are also called Euclidean geometry codes, while the codes based on $PG_s(m,q)$
are called projective geometry codes.
The codes over $GF(p)$ thus defined, where $q=p^t$,
correspond to subfield subcodes of generalized Reed-Muller codes
\cite[Chapter 5]{AK},  \cite{GD2}.
The binary Euclidean geometry code
being the null space of the incidence matrix of $AG_s(m,2)$
is equivalent
to the Reed-Muller code $R(m-s,m)$ of length $2^m$ and order $m-s$.
It is well known that the finite geometry codes admit majority logic decoding
 \cite{Reed}, \cite{W}, \cite{GD}.

In \cite{JT09}, Jungnickel and Tonchev
 used polarities in projective geometry
to find a class of designs which have the same parameters
and share some other properties with a projective geometry design $PG_s(2s,q)$, $s\ge 2$,
but are not isomorphic to $PG_s(2s,q)$.
We refer to these designs as {\bf polarity designs}. In the cases when $q=p$ is a prime,
the $p$-rank of the incidence matrix of a polarity design $D$ is equal to that of $PG_s(2s,p)$,
hence the polarity designs provide an infinite class of counter-examples to Hamada's 
conjecture \cite{H1}, \cite{H2}.

In \cite{CT}, Clark and Tonchev proved that the code being the null space of
the incidence matrix of a polarity design can correct by majority-logic decoding the
same number of errors as the projective geometry code based on $PG_s(2s,q)$.
In the binary case ($q=2$), the minimum distance of the code of the polarity design
obtained from $PG(2s,2)$, is
equal to $2^{s+1}$, and the majority-logic algorithm from \cite{CT} corrects 
all errors guaranteed by the
minimum distance. The extended code of the binary code spanned by
the blocks of a polarity design obtained from $PG(2s,2)$ is a self-dual
binary code of the same length, dimension and minimum distance, and correcting by
majority-logic the same number of errors as the Reed-Muller code
$R(s,2s+1)$ of length $2^{2s+1}$ and order $s$.

In the smallest binary case, $s=2$, the extended code of the polarity design
obtained from $PG(4,2)$, 
is a doubly-even self-dual $[32,16,8]$ code, which not only has the same parameters
and corrects by majority-logic decoding the same number of errors
as the 2nd order Reed-Muller code $R(2,5)$, but also has the same weight distribution
as $R(2,5)$. This phenomenon is easily explained by the fact that both codes are
extremal doubly-even self-dual codes, hence are forced to have the same weight distribution
\cite{NRS} (actually, in this case there are five inequivalent extremal
doubly-even self-dual
$[32,16,8]$ codes \cite{CP}.).

It is the aim of this note to report the computation of the weight distribution
of the extended code of the polarity design in the next case $s=3$,
i.e.\ the polarity design obtained
from $PG(6,2)$, and to demonstrate that this doubly-even  self-dual $[128,64,16]$ code
has the same weight distribution as the 3rd order Reed-Muller code $R(3,7)$. 

One of the authors, Vladimir Tonchev, conjectures that the extended code
of the polarity design obtained from $PG(2s,2)$ has the same weight distribution as
the Reed-Muller code $R(s,2s+1)$ for every $s\ge 2$.

\section{Computing the weight distribution}

The polarity design $D$ obtained from $PG(6,2)$ \cite{JT09}
is a 2-$(v,k,\lambda)$ design with parameters 
\[ v=127, \ k=15, \ \lambda=155, \]
that is, $D$ has the same parameters as the projective geometry
design $PG_{3}(6,2)$ having as blocks the 3-dimensional subspaces
of the 6-dimensional projective geometry $PG(6,2)$ over the field of order 2.
In addition, ${D}$ has the same block intersections as  $PG_{3}(6,2)$,
namely 1, 3 and 7, and its incidence matrix has the same 2-rank 64 as
 $PG_{3}(6,2)$ \cite{JT09},
hence provides a counter-example to Hamada's conjecture \cite{H1}, \cite{H2}.

The parameters and the block intersection numbers of ${D}$ imply that
the binary linear code $C$ spanned by the block by point incidence matrix of ${D}$
has minimum distance not exceeding 15, and its extended code $C^{*}$ is a doubly-even
self-dual $[128,64]$ code of minimum distance $d\le 16$.
It follows from the results from \cite{CJT} and \cite{CT} that $d=16$ and the
 code $C^{*}$ admits
majority-logic decoding that corrects up to 7 errors, that is, the same number of errors
as the doubly-even self-dual $[128,64,16]$ 3rd order Reed-Muller code $R(3,7)$.

We will show that $C^{*}$ has the same weight distribution as $R(3,7)$.
The weight distribution of the Reed-Muller code $R(3,7)$ was computed by
Sugino, Ienaga,  Tokura and  Kasami \cite{SITK} (see 
{\it The On-line Encyclopedia of Integer Sequences} \cite{OEIS}, sequence A110845),
and is listed in Table \ref{RM}.

\begin{table}[htpb!]
\begin{center}
\begin{tabular}{|r|l|}
\hline
$i$ & $A _i$ \\
\hline
0 & 1\\

16 & 94488 \\

20 & 0\\

24 & 74078592\\

28 & 3128434688\\

32 & 312335197020\\

36 & 18125860315136\\

40 & 552366841342848\\

44 & 9491208609103872\\

48 & 94117043084875944\\

52 & 549823502398291968\\

56 & 1920604779257215744\\

60 & 4051966906789380096\\

64 & 5193595576952890822\\

68 & 4051966906789380096\\

72 & 1920604779257215744\\

76 & 549823502398291968\\

80 & 94117043084875944 \\

84 & 9491208609103872 \\

88 & 552366841342848\\

92 & 18125860315136 \\

96 & 312335197020 \\

100 & 3128434688 \\

104 & 74078592 \\

108 & 0\\

112 & 94488 \\

128 & 1\\
\hline

\end{tabular}\caption{Weight distribution of $R(3,7)$}
\label{RM}
\end{center}
\end{table}

Since the code dimension 64 is significant, to facilitate the computation 
of the weight distribution of $C^{*}$, we
employ known properties of weight enumerators of binary doubly-even self-dual codes.
Since $C^{*}$ is a doubly-even self-dual  $[128,64,16]$ code,
by the Gleason theorem (cf.~\cite[Section 2]{MaS}, \cite{NRS}),
the weight distribution $\{ A_i \}_{i=0}^{128}$
of  $C^{*}$ can be determined completely by the values of $A_{16}$ and $A_{20}$.
More specifically, the two-variable weight enumerator can be written as
\[  
\sum_{i=0}^{128} A_{i}x^{128-i}y^i 
=\sum_{j=0}^5 b_{j}(x^8 +14x^{4}y^4 +y^8)^{16-3j}
(x^4y^4(x^4-y^4)^4)^j, 
\]
where
\begin{align*}
&
b_0=1,
b_1=-224,
b_2=16336,
b_3=-430656,
\\&
b_4=A_{16}+3196776 \text{ and }
b_5=A_{20}-40A_{16}-2696256.
\end{align*}
Consequently, the weight enumerator 
\[ W(x) = \sum_{i=0}^{128} A_{i}x^i \]
can be written as
\begin{align*}
 W(x) = & 1+A_{16}x^{16} + A_{20}x^{20} \\
&+ (13228320 + 644A_{16} - 6A_{20})x^{24}\\
&+(2940970496 + 1984A_{16} - 89A_{20})x^{28} \\
&+  (320411086380 - 85470A_{16} + 1500A_{20})x^{32} + \cdots
\end{align*}
 (cf. \cite[Section 2]{MaS}).

A $64 \times 128$ generator matrix $G$ of the extended code $C^{*}$ 
was computed following the construction of polarity designs from \cite{JT09},
and is available on-line at

\begin{verbatim}
http://www.math.mtu.edu/~tonchev/bordered_genmatrix.txt
\end{verbatim}

Using Magma \cite{Magma},  it took a few minutes to compute on a PC that 
\begin{equation}
\label{a16}
 A_{16}=94488, 
\end{equation}
and about an hour\footnote{It is an interesting question if one can prove that $A_{20}=0$ geometrically,
using the construction from \cite{JT09}.}
 to compute 
\begin{equation}
\label{a20}
 A_{20}=0.
\end{equation}

Since the values $A_{16}$ (see (\ref{a16})) and $A_{20}$ (see (\ref{a20})) are the same as
the corresponding values for the self-dual $[128,64,16]$  Reed-Muller code $R(3,7)$
(cf.~Table \ref{RM}), we have the following.

\begin{theorem}\label{thm}
The weight distribution of the extended $[128,64,16]$ code $C^{*}$ of the the code
$C$ spanned by the incidence vectors of the blocks of the polarity design $D$
obtained from $PG(6,2)$, is identical with the weight distribution 
of the 3rd order Reed-Muller code $R(3,7)$.
\end{theorem}
 
Using Magma, it took 90 seconds to compute the full automorphism group
$\operatorname{Aut}(C^{*})$
of $C^{*}$. Since $C^{*}$ is spanned by the set of minimum weight vectors which form
the block by point incidence matrix of a 3-$(128,16,155)$ design $D^{*}$
 \cite{CJT}, \cite{CT}, the full automorphism group of $C^{*}$ coincides with that of
$D^{*}$, and is of order
\begin{equation}
\label{autC*}
|\operatorname{Aut}(C^{*})|=165140150353920 =  2^{28} \cdot 3^4 \cdot 5 \cdot 7^2 \cdot 31.
\end{equation}
It follows from results of \cite{JT09} that the polarity 2-$(127,15,155)$ design $D$,
and hence, the related $[127,64]$ code $C$, is invariant under the collineation
 group of the affine space $AG(6,2)$, being  of order
\[ 2^6(2^6 -1)(2^6 -2)(2^6-2^2)(2^6-2^3)(2^6-2^4)(2^6-2^5)= 2^{21}\cdot 3^4
 \cdot 5 \cdot 2^7 \cdot 31. \]
 Thus, the full automorphism group $\operatorname{Aut}(C^{*})$ of $C^{*}$  
extends the automorphism group of $C$ by a factor of $2^7 =128$,
and acts transitively on the set of 128 code coordinates.

It is known that the full automorphism group $\operatorname{Aut}(R(3,7))$
of the Reed-Muller code $R(3,7)$ is equivalent to the group of collineations of $AG(7,2)$,
and is of order
\begin{equation}
\label{autRM}
|\operatorname{Aut}(R(3,7))|=20972799094947840 = 2^{28} \cdot 3^4 \cdot 5 \cdot 7^2 \cdot 31 \cdot 127.
\end{equation}

\begin{remark}
{\rm
The parameters of the extended self-dual
code, obtained from the polarity design in $PG(8,2)$,
are $[512,256,32]$. It seems computationally infeasible to find the weight
distribution of such a code by computer, 
even with the help of Gleason's theorem, 
due to the very large code dimension. 
 Thus, any proof of the conjecture formulated
in the last paragraph of Introduction,
has to be based on geometric or other theoretical considerations.
}
\end{remark}

\section{Self-dual neighbors}

In this section, we investigate self-dual neighbors of 
the $3$rd order Reed-Muller code $R(3,7)$
and the extended $[128,64,16]$ code $C^{*}$ given in Theorem \ref{thm}.
Two self-dual codes $C$ and $C'$ of length $n$
are said to be {neighbors} if $\dim(C \cap C')=n/2-1$. 
We give some observations from \cite{CHK64}
on self-dual codes constructed by neighbors.
Let $C$ be a self-dual $[n,n/2,d]$ code.
Let $M$ be a matrix  whose rows are the codewords of weight $d$ in $C$.
Suppose that there is a vector $x\in GF(2)^n$ such that
\begin{equation}\label{Eq:N}
M x^T = \allone^T,
\end{equation}
where 
$x^T$ denotes the transpose of $x$ and 
$\allone$ is the all-one vector.
Set $C_0=\langle x\rangle^\perp\cap C$, where $\langle x\rangle$ denotes
the code generated by $x$.
Then $C_0$ is a subcode of index $2$ in $C$.
If the weight of $x$ is even, then
we have the two self-dual neighbors $\langle C_0,x\rangle$ and $\langle
C_0,x+y\rangle$ of $C$
for some $y\in C\setminus C_0$,
which do not have any codeword of weight $d$ in $C$,
where $\langle C,x\rangle = C \cup (x+C)$. 
When $C$ has a self-dual $[n,n/2,d']$ neighbor $C'$ with $d'\ge d+2$,
(\ref{Eq:N}) has a solution $x$ and we can obtain $C'$ in this way.
If $\rank M < \rank ( M\ \allone^T)$,
then $C$ has no self-dual $[n,n/2,d']$ neighbor $C'$ with $d'\ge d+2$.
Using Magma, we verified that
\[
(\rank M,\rank ( M\ \allone^T))=(64, 65)
\]
for the 3rd order Reed-Muller code $R(3,7)$
and the extended $[128,64,16]$ code $C^{*}$ given in Theorem \ref{thm}.
Therefore, we have the following:

\begin{proposition}
The $3$rd order Reed-Muller code $R(3,7)$
and the extended $[128,64,16]$ code $C^{*}$ given in Theorem \ref{thm}
have no self-dual $[128,64,d]$ neighbor
with $d \in\{20,24\}$.
\end{proposition}

This means that the above two doubly-even self-dual codes
of length $128$ have no extremal doubly-even self-dual neighbor
of that length.

\section{ Acknowledgments}

The authors thank the anonymous reviewers for the useful comments
and suggestions that lead to some improvements of the text.

\end{document}